\documentclass[12pt,leqno,draft]{article}

\def\diam{\mathop{\rm diam}}
\def\dist{\mathop{\rm dist}}

\begin{document}

\title{Some remarks about Cauchy integrals}

\author{Stephen William Semmes  \\
	Rice University		\\
	Houston, Texas}

\date{}

\maketitle

	A basic theme in the wonderful books and surveys of Stein,
Weiss, and Zygmund is that Hilbert transforms, Poisson kernels, heat
kernels, and related objects are quite interesting and fundamental.  I
certainly like this point of view.  There is a variety of ways in
which things can be interesting or fundamental, of course.

	In the last several years there have been striking
developments connected to Cauchy integrals, and in this regard I would
like to mention the names of Pertti Mattila, Mark Melnikov, and Joan
Verdera in particular.  I think many of us are familiar with the
remarkable new ideas involving Menger curvature, and indeed a lot of
work using this has been done by a lot of people, and continues to be
done.  Let us also recall some matters related to \emph{symmetric
measures}.

	Let $\mu$ be a nonnegative Borel measure on the complex plane
${\bf C}$, which is finite on bounded sets.  Following Mattila, $\mu$
is said to be symmetric if for each point $a$ in the support of $\mu$
and each positive real number $r$ we have that the integral of $z - a$
over the open ball with center $a$ and radius $r$ with respect to the
measure $d\mu(z)$ is equal to $0$.  One might think of this as a kind
of flatness condition related to the existence of principal values of
Cauchy integrals.

	If $\mu$ is equal to a constant times $1$-dimensional Lebesgue
measure on a line, then $\mu$ is a symmetric measure.  For that
matter, $2$-dimensional Lebesgue measure on ${\bf C}$ is symmetric,
and there are other possibilities.  Mattila discusses this, and shows
that a symmetric measure which satisfies some additional conditions is
equal to a constant times $1$-dimensional Lebesgue measure on a line.

	Mattila uses this to show that existence almost everywhere of
principle values of a measure implies rectifiability properties of the
measure.  Mattila and Preiss have considered similar questions in
general dimensions, where the geometry is more complecated.  Mattila's
student Petri Huovinen has explored analogous matters in the plane
with more tricky kernels than the Cauchy kernel.

	Another kind of $m$-dimensional symmetry condition for a
nonnegative Borel measure $\mu$ on ${\bf R}^n$, which is finite on
bounded sets, asks that for each point $a$ in the support of $\mu$ and
for each radius $r > 0$ the $\mu$ measure of the open ball with center
$a$ and radius $r$ is equal to a constant $c$, depending only on
$\mu$, times $r^m$.  This condition holds for constant multiples of
$m$-dimensional Lebesgue measure on $m$-dimensional planes in ${\bf
R}^n$.  The converse is known in some cases, and non-flat examples are
also known in some cases.  These types of measures have been studied
extensively by Preiss, partly in collaboration with Kirchheim and with
Kowalski, and it seems fair to say that many mysteries remain.

	In general, it seems to me that there are a lot of very
interesting questions involving geometry of sets, measures, currents,
and varifolds and quantities such as Cauchy integrals, measurements of
symmetry like those considered by Huovinen, Mattila, and Preiss, and
densities ratios.  This may entail rather exact conditions and special
geometric structures, or approximate versions and some kind of
regularity.  In the latter case, instead of asking that some quantity
vanish exactly, one might look at situations where it satisfies a
bound like $O(r^\alpha)$ for some positive real number $\alpha$, and
where $r > 0$ is a radius or similar parameter.

	I like very much a paper by Verdera on $T(1)$ theorems for
Cauchy integrals, which uses Menger curvature ideas.  It seems to me
that it could be a starting point for a new kind of operator theory
for certain kinds of operators.  There is a lot of room for
development for new kinds of structure of linear operators.

	Another reasonably-specific area with a lot of possibilities
is to try to combine Menger curvature ideas with the rotation method.
At first they may not seem to fit together too easily.  However, I
would not be too surprised if some interesting things could come up in
this manner.

\section{Integrals of curvature on curves and surfaces}\footnotetext{A 
presentation based on this section was made at the AMS Special Session
``Surface Geometry and Shape Perception'' (Hoboken, 2001).}
\setcounter{equation}{0}

	In this section we discuss some topics that came up in
Chapters 2 and 3 of Part III of \cite{DS2}.  These involve relations
between derivatives of Cauchy integrals on curves and surfaces and
curvatures of the curves and surfaces.  In ${\bf R}^n$ for $n > 2$,
``Cauchy integrals'' can be based on generalizations of complex
analysis using quarternions or Clifford algebras (as in \cite{BDS}).
Part of the point here is to bring out the basic features and types of
computations in a simple way, if not finer aspects which can also be
considered.

	Let us consider first curves in the plane ${\bf R}^2$.  We
shall identify ${\bf R}^2$ with the set ${\bf C}$ of complex numbers.

	Let $\Gamma$ be some kind of curve in ${\bf C}$, or perhaps
union of pieces of curves.  For each $z \in {\bf C} \backslash \Gamma$,
we have the contour integral
\begin{equation}
\label{int_{Gamma} frac{1}{(z-zeta)^2} d zeta}
	\int_{\Gamma} \frac{1}{(z-\zeta)^2} \, d\zeta
\end{equation}
as from complex analysis.  More precisely, ``$d\zeta$'' is the element
of integration such that if $\gamma$ is an arc in ${\bf C}$ from a
point $a$ to another point $b$, then
\begin{equation}
	\int_{\gamma} d\zeta = b - a.
\end{equation}
This works no matter how $\gamma$ goes from $a$ to $b$.  This is different
from integrating arclength, for which the element of integration is often
written $|d\zeta|$.  For this we have that
\begin{equation}
	\int_{\gamma} |d\zeta| = {\rm length}(\gamma),
\end{equation}
and this very much depends on the way that $\gamma$ goes from $a$ to
$b$.

	If $\Gamma$ a union of closed curves, then
\begin{equation}
\label{int_{Gamma} frac{1}{(z-zeta)^2} d zeta = 0}
	\int_{\Gamma} \frac{1}{(z-\zeta)^2} \, d\zeta = 0.
\end{equation}
This is a standard formula from complex analysis (an instance of
``Cauchy formulae''), and one can look at it in the following manner.
As a function of $\zeta$, $1/(z-\zeta)^2$ is the complex derivative in
$\zeta$ of $1/(z-\zeta)$,
\begin{equation}
	\frac{d}{d\zeta} \biggl(\frac{1}{z-\zeta}\biggr)
				= \frac{1}{(z-\zeta)^2}.
\end{equation}
If $\gamma$ is a curve from points $a$ to $b$ again, which does not pass
through $z$, then
\begin{equation}
\label{int_{gamma} frac{1}{(z-zeta)^2} d zeta = frac{1}{z-b} - frac{1}{z-a}}
	\int_{\gamma} \frac{1}{(z-\zeta)^2} \, d\zeta
		= \frac{1}{z-b} - \frac{1}{z-a}.
\end{equation}
In particular, one gets $0$ for closed curves (since that corresponds to
having $a = b$).

	As a variation of these matters, if $\Gamma$ is a line, then
\begin{equation}
\label{int_{Gamma} frac{1}{(z-zeta)^2} d zeta = 0, 2}
	\int_{\Gamma} \frac{1}{(z-\zeta)^2} \, d\zeta = 0
\end{equation}
again.  This can be derived from (\ref{int_{gamma} frac{1}{(z-zeta)^2}
d zeta = frac{1}{z-b} - frac{1}{z-a}}) (and can be looked at in terms
of ordinary calculus, without complex analysis).  There is enough
decay in the integral so that there is no problem with using the whole
line.

	What would happen with these formulae if we replaced the
complex element of integration $d\zeta$ with the arclength element of
integration $|d\zeta|$?  In general we would not have
(\ref{int_{Gamma} frac{1}{(z-zeta)^2} d zeta = 0}) for unions of
closed curves, or (\ref{int_{gamma} frac{1}{(z-zeta)^2} d zeta =
frac{1}{z-b} - frac{1}{z-a}}) for a curve $\gamma$ from $a$ to $b$.
However, we would still have (\ref{int_{Gamma} frac{1}{(z-zeta)^2} d
zeta = 0, 2}) for a line, because in this case $d\zeta$ would be a
constant times $|d\zeta|$.

	Let us be a bit more general and consider an element of
integration $d\alpha(\zeta)$ which is positive, like the arclength
element $|d\zeta|$, but which is allowed to have variable density.
Let us look at an integral of the form
\begin{equation}
\label{int_{Gamma} frac{1}{(z-zeta)^2} dalpha(zeta)}
	\int_{\Gamma} \frac{1}{(z-\zeta)^2} \, d\alpha(\zeta).
\end{equation}
This integral can be viewed as a kind of measurement of \emph{curvature}
of $\Gamma$ (which also takes into account the variability of the density
in $d\alpha(\zeta)$).

	If we put absolute values inside the integral, then the result
would be roughly $\dist(z,\Gamma)^{-1}$,
\begin{equation}
\label{int_{Gamma} frac{1}{|z-zeta|^2} dalpha(zeta) approx dist(z,Gamma)^{-1}}
	\int_{\Gamma} \frac{1}{|z-\zeta|^2} \, d\alpha(\zeta)
			\approx \dist(z,\Gamma)^{-1}
\end{equation}
under suitable conditions on $\Gamma$.  For instance, if $\Gamma$ is a
line, then the left side of (\ref{int_{Gamma} frac{1}{|z-zeta|^2}
dalpha(zeta) approx dist(z,Gamma)^{-1}}) is equal to a positive
constant times the right side of (\ref{int_{Gamma} frac{1}{|z-zeta|^2}
dalpha(zeta) approx dist(z,Gamma)^{-1}}).  

	The curvature of a curve is defined in terms of the derivative
of the unit normal vector along the curve, or, what is essentially the
same here, the derivative of the unit tangent vector.  The unit tangent
vector gives exactly what is missing from $|d\zeta|$ to get $d\zeta$,
if we write the unit tangent vector as a complex number.  (One should
also follow the tangent in the orientation of the curve.)

	If the curve is a line or a line segment, then the tangent
is constant, which one can pull in and out of the integral.  In general
one can view (\ref{int_{Gamma} frac{1}{(z-zeta)^2} dalpha(zeta)}) as
a measurement of the variability of the unit tangent vectors, and of
the variability of the positive density involved in $d\alpha(\zeta)$.

	Let us look at some simple examples.  Suppose first that $\Gamma$
is a straight line segment from a point $a \in {\bf C}$ to another point
$b$, $a \ne b$.  Then $|d\zeta|$ is a constant multiple of $d\zeta$, and
\begin{equation}
	\int_{\Gamma} \frac{1}{(z-\zeta)^2} \, |d\zeta|
	   = ({\rm constant}) \cdot \Bigl(\frac{1}{z-b} - \frac{1}{z-a}\Bigr).
\end{equation}
In this case the ordinary curvature is $0$, except that one can say
that there are contributions at the endpoints, like Direc delta
functions, which are reflected in right side.  If $z$ gets close to
$\Gamma$, but does not get close to the endpoints $a$, $b$ of
$\Gamma$, then the right side stays bounded and behaves nicely.  This
is ``small'' in comparison with $\dist(z,\Gamma)^{-1}$.  Near $a$ or
$b$, we get something which is indeed like $|z-a|^{-1}$ or
$|z-b|^{-1}$.

	As another example, suppose that we have a third point $p \in
{\bf C}$, where $p$ does not lie in the line segment between $a$ and
$b$ (and is not equal to $a$ or $b$).  Consider the curve $\Gamma$
which goes from $a$ to $p$ along the line segment between them, and
then goes from $p$ to $b$ along the line segment between them.  Again
$|d\zeta|$ is a constant multiple of $d\zeta$.  Now we have
\begin{equation}
	\int_{\Gamma} \frac{1}{(z-\zeta)^2} \, |d\zeta|
	   = c_1 \Bigl(\frac{1}{z-p} - \frac{1}{z-a}\Bigr)
		+ c_2 \Bigl(\frac{1}{z-b} - \frac{1}{z-p}\Bigr),
\end{equation}
where $c_1$ and $c_2$ are constants which are not equal to each other.
This is like the previous case, except that the right side behaves like
a constant times $|z-p|^{-1}$ near $p$ (and remains bounded away from
$a$, $b$, $p$).  This reflects the presence of another Dirac delta
function for the curvature, at $p$.  If the curve flattens out, so that
the angle between the two segments is close to $\pi$, then the coefficient
$c_1 - c_2$ of the $(z-p)^{-1}$ term becomes small.

	Now suppose that $\Gamma$ is the unit circle in ${\bf C}$, centered
around the origin.  In this case $|d\zeta|$ is the same as $d\zeta/\zeta$,
except for a constant factor, and we consider
\begin{equation}
	\int_{\Gamma} \frac{1}{(z-\zeta)^2} \, \frac{d\zeta}{\zeta}.
\end{equation}
If $z = 0$, then one can check that this integral is $0$.  For $z \ne 0$,
let us rewrite the integral as
\begin{equation}
	\int_{\Gamma} \frac{1}{(z-\zeta)^2} \, 
			\Bigl(\frac{1}{\zeta} - \frac{1}{z}\Bigr) \, d\zeta
	+ \frac{1}{z} \int_{\Gamma} \frac{1}{(z-\zeta)^2} \, d\zeta.
\end{equation}
The second integral is $0$ for all $z \in {\bf C} \backslash \Gamma$,
as in the earlier discussion.  The first integral is equal to
\begin{equation}
   \int_{\Gamma} \frac{1}{(z-\zeta)^2} \, \frac{(z-\zeta)}{z \zeta} \, d\zeta
= \frac{1}{z} \int_{\Gamma} \frac{1}{(z-\zeta)} \, \frac{1}{\zeta} \, d\zeta.
\end{equation}
On the other hand,
\begin{equation}
	\frac{1}{(z-\zeta)} \, \frac{1}{\zeta} 
	    = \frac{1}{z} \Bigl(\frac{1}{z-\zeta} + \frac{1}{\zeta}\Bigr),
\end{equation}
and so we obtain 
\begin{equation}
   \frac{1}{z^2} \int_{\Gamma} \Bigl(\frac{1}{z-\zeta} + \frac{1}{\zeta}\Bigr)
								 \, d\zeta.
\end{equation}
For $|z| > 1$ we have that
\begin{equation}
	\int_{\Gamma} \frac{1}{z-\zeta} \, d\zeta = 0,
\end{equation}
and thus we get a constant times $1/z^2$ above.  If $|z| < 1$, then
\begin{equation}
	\int_{\Gamma} \frac{1}{z-\zeta} \, d\zeta 
		= - \int_{\Gamma} \frac{1}{\zeta} \, d\zeta,
\end{equation}
and the earlier expression is equal to $0$.

	For another example, fix a point $q \in {\bf C}$, and suppose that
$\Gamma$ consists of a finite number of rays emanating from $q$.  On
each ray, we assume that we have an element of integration $d\alpha(\zeta)$
which is a positive constant times the arclength element.  

	If $R$ is one of these rays, then 
\begin{equation}
	\int_R \frac{1}{(z-\zeta)^2} \, d\alpha(\zeta)
		= ({\rm constant}) \, \frac{1}{z-q}.
\end{equation}
This constant takes into account both the direction of the ray and the
density factor in $d\alpha(\zeta)$ on $R$.

	If we now sum over the rays, we still get 
\begin{equation}
	\int_{\Gamma} \frac{1}{(z-\zeta)^2} \, d\alpha(\zeta)
		= ({\rm constant}) \, \frac{1}{z-q};
\end{equation}
however, this constant can be $0$.  This happens if $\Gamma$ is a
union of lines through $q$, with constant density on each line, and it
also happens more generally, when the directions of the rays satisfy a
suitable balancing condition, depending also on the density factors
for the individual rays.  This can happen with $3$ rays, for instance.

	When the constant is $0$, $\Gamma$ (with these choices of
density factors) has ``curvature $0$'', even if this is somewhat
complicated, because of the singularity at $q$.  This is a special
case of the situation treated in \cite{AA}.

	In general, ``weak'' or integrated curvature is defined using
suitable test functions on ${\bf R}^2$ with values in ${\bf R}^2$ (or
on ${\bf R}^n$ with values in ${\bf R}^n$), as in \cite{AA}.  For $n =
2$ one can reformulate this in terms of complex-valued functions on
${\bf C}$, and complex-analyticity gives rise to simpler formulas.
The link between this kind of story with Cauchy integrals and the weak
notion of curvature for varifolds as in \cite{AA} was suggested by Bob
Hardt.

	For more information on these topics, see Chapter 2 of Part
III of \cite{DS2}.  In \cite{DS2} there are further issues which are not
needed in various settings.
	
	Now let us look at similar matters in ${\bf R}^n$, $n > 2$,
and $(n-1)$-dimensional surfaces there.  Ordinary complex analysis is
no longer available, but there are substitutes, in terms of
quarternions (in low dimensions) and Clifford algebras.  For the sake
of definiteness let us focus on the latter.

	Let $n$ be a positive integer.  The \emph{Clifford algebra}
$\mathcal{C}(n)$ has $n$ generators $e_1, e_2, \ldots, e_n$ which
satisfy the following relations:
\begin{eqnarray}
	e_j \, e_k & = & - e_k \, e_j \quad\hbox{when } j \ne k;	\\
	e_j^2 & = & -1	       \qquad\hbox{ for all } j.	\nonumber
\end{eqnarray}
Here $1$ denotes the identity element in the algebra.  These are the
only relations.  More precisely, one can think of $\mathcal{C}(n)$ first
as a real vector space of dimension $2^n$, in which one has a basis
consisting of all products of $e_j$'s of the form
\begin{equation}
	e_{j_1} \, e_{j_2} \cdots e_{j_\ell},
\end{equation}
where $j_1 < j_2 < \cdots < j_\ell$, and $\ell$ is allowed to range from
$0$ to $n$, inclusively.  When $\ell = 0$ this is interpreted as
giving the identity element $1$.  If $\beta, \gamma \in
\mathcal{C}(n)$, then $\beta$ and $\gamma$ are given by linear
combinations of these basis elements, and it is easy to define the
product $\beta \, \gamma$ using the relations above and standard rules
(associativity and distributivity).

	If $n = 1$, then the result is isomorphic to the complex
numbers in a natural way, and if $n = 2$, the result is isomorphic to
the quarternions.  Note that $\mathcal{C}(n)$ contains ${\bf R}$
in a natural way, as multiples of the identity element.

	A basic feature of the Clifford algebra $\mathcal{C}(n)$ is
that if $\beta \in \mathcal{C}(n)$ is in the linear span of $e_1, e_2,
\ldots, e_n$ (without taking products of the $e_j$'s), then $\beta$
can be inverted in the algebra if and only if $\beta \ne 0$.  More
precisely, if
\begin{equation}
	\beta = \sum_{j=1}^n \beta_j \, e_j,
\end{equation}
where each $\beta_j$ is a real number, then
\begin{equation}
	\beta^2 = - \sum_{j=1}^n |\beta_j|^2.
\end{equation}
If $\beta \ne 0$, then the right side is a nonzero real number,
and $-(\sum_{j=1}^n |\beta_j|^2)^{-1} \beta$ is the multiplicative
inverse of $\beta$.

	More generally, if $\beta$ is in the linear span of
$1$ and $e_1, e_2, \ldots, e_n$, so that
\begin{equation}
	\beta = \beta_0 + \sum_{j=1}^n \beta_j \, e_j,
\end{equation}
where $\beta_0, \beta_1, \ldots, \beta_n$, then we set
\begin{equation}
	\beta^* = \beta_0 - \sum_{j=1}^n \beta_j \, e_j.
\end{equation}
This is analogous to complex conjugation of complex numbers, and
we have that
\begin{equation}
	\beta \, \beta^* = \beta^* \, \beta = \sum_{j=0}^n |\beta_j|^2.
\end{equation}
If $\beta \ne 0$, then $(\sum_{j=0}^n |\beta_j|^2)^{-1} \beta^*$ is
the multiplicative inverse of $\beta$, just as in the case of complex
numbers.

	When $n > 2$, nonzero elements of $\mathcal{C}(n)$ may not be
invertible.  For real and complex numbers and quarternions it is true
that nonzero elements are invertible.  The preceding observations are
substitutes for this which are often sufficient.

	Now let us turn to \emph{Clifford analysis}, which is an
analogue of complex numbers in higher dimensions using Clifford
algebras.  (See \cite{BDS} for more information.)

	Suppose that $f$ is a function on ${\bf R}^n$, or some
subdomain of ${\bf R}^n$, which takes values in $\mathcal{C}(n)$.  We
assume that $f$ is smooth enough for the formulas that follow (with
the amount of smoothness perhaps depending on the circumstances).
Define a differential operator $\mathcal{D}$ by
\begin{equation}
	\mathcal{D} f = \sum_{j=1}^n e_j \, \frac{\partial}{\partial x_j} f.
\end{equation}

	Actually, there are some natural variants of this to also
consider.  This is the ``left'' version of the operator; there is also
a ``right'' version, in which the $e_j$'s are moved to the right side of
the derivatives of $f$.  This makes a difference, because the Clifford
algebra is not commutative, but the ``right'' version enjoys the same
kind of properties as the ``left'' version.  (Sometimes one uses the
two at the same time, as in certain integral formulas, in which the
two operators are acting on separate functions which are then part of
the same expression.)  

	As another alternative, one can use the Clifford algebra
$\mathcal{C}(n-1)$ for Clifford analysis on ${\bf R}^n$, with one
direction in ${\bf R}^n$ associated to the multiplicative identity
element $1$, and the remaining $n-1$ directions associated to $e_1,
e_2, \ldots, e_{n-1}$.  There is an operator analogous to
$\mathcal{D}$, and properties similar to the ones that we are about to
describe (with adjustments analogous to the conjugation operation
$\beta \mapsto \beta^*$).

	For the sake of definiteness, let us stick to the version that
we have.  A function $f$ as above is said to be \emph{Clifford analytic}
if
\begin{equation}
	\mathcal{D} f = 0
\end{equation}
(on the domain of $f$).

	Clifford analytic functions have a lot of features analogous
to those of complex analytic functions, including integral formulas.
There is a natural version of a \emph{Cauchy kernel}, which is given
by
\begin{equation}
	\mathcal{E}(x-y) = \frac{\sum_{j=1}^n (x_j - y_j) \, e_j}{|x-y|^n}.
\end{equation}
This function is Clifford analytic in $x$ and $y$ away from $x = y$,
and it has a ``fundamental singularity'' at $x = y$, just as $1/(z-w)$
has in the complex case.

	One can calculate these properties directly, and one can also
look at them in the following way.  A basic indentity involving
$\mathcal{D}$ is 
\begin{equation}
	\mathcal{D}^2 = - \Delta,
\end{equation}
where $\Delta$ denotes the Laplacian, $\Delta = \sum_{j=1}^n \partial^2 /
\partial x_j^2$.  The kernel $\mathcal{E}(x)$ is a constant multiple of
\begin{eqnarray}
	&& \mathcal{D}(|x|^{n-2}) \quad\hbox{when } n > 2,		\\
	&& \mathcal{D}(\log |x|)  \quad\hbox{when } n = 2.	\nonumber
\end{eqnarray}
For instance, the Clifford analyticity of $\mathcal{E}(x)$ for $x \ne 0$
follows from the harmonicity of $|x|^{n-2}$, $\log |x|$ for $x \ne 0$
(when $n > 2$, $n = 2$, respectively).

	Analogous to (\ref{int_{Gamma} frac{1}{(z-zeta)^2} d zeta}), let
us consider integrals of the form
\begin{equation}
\label{int_{Gamma} frac{partial}{partial x_m} mathcal{E}(x-y) N(y) dy}
	\int_{\Gamma} \frac{\partial}{\partial x_m} \mathcal{E}(x-y) \,
 		     N(y) \, dy,   \quad x \in {\bf R}^n \backslash \Gamma,
\end{equation}
where $\Gamma$ is some kind of $(n-1)$-dimensional surface in ${\bf R}^n$,
or union of pieces of surfaces, 
\begin{equation}
	N(y) = \sum_{j=1}^n N_j(y) \, e_j
\end{equation}
is the unit normal to $\Gamma$ (using some choice of orientation for
$\Gamma$), turned into an element of $\mathcal{C}(n)$ using the
$e_j$'s in this way, and $dy$ denotes the usual element of surface
integration on $\Gamma$.  Thus $N(y) \, dy$ is a
Clifford-algebra-valued element of integration on $\Gamma$ which is
analogous to $d\zeta$ for complex contour integrals, as in
(\ref{int_{Gamma} frac{1}{(z-zeta)^2} d zeta}).  A version of the
Cauchy integral formula implies that
\begin{equation}
	\int_{\Gamma} \mathcal{E}(x-y) \, N(y) \, dy
\end{equation}
is locally constant on ${\bf R}^n \backslash \Gamma$ when $\Gamma$ is
a ``closed surface'' in ${\bf R}^n$, i.e., the boundary of some
bounded domain (which is reasonably nice).  In fact, this integral is
a nonzero constant inside the domain, and it is zero outside the
domain.  At any rate, the differentiated integral (\ref{int_{Gamma}
frac{partial}{partial x_m} mathcal{E}(x-y) N(y) dy}) is then $0$ for
all $x \in {\bf R}^n \backslash \Gamma$, in analogy with
(\ref{int_{Gamma} frac{1}{(z-zeta)^2} d zeta = 0}).

	Now suppose that we have a positive element of integration
$d\alpha(y)$ on $\Gamma$, which is the usual element of surface
integration $dy$ together with a positive density which is allowed
to be variable.  Consider integrals of the form
\begin{equation}
\label{int_{Gamma} frac{partial}{partial x_m} mathcal{E}(x-y) d alpha(y)}
	\int_{\Gamma} \frac{\partial}{\partial x_m} \mathcal{E}(x-y) 
 		     \, d\alpha(y),   \quad x \in {\bf R}^n \backslash \Gamma.
\end{equation}
This again can be viewed in terms of integrations of curvatures of
$\Gamma$ (also incorporating the variability of the density in
$d\alpha(y)$).  In a ``flat'' situation, as when $\Gamma$ is an
$(n-1)$-dimensional plane, or a piece of one, $N(y)$ is constant, and
if $d\alpha(y)$ is replaced with a constant times $dy$, then we can
reduce to (\ref{int_{Gamma} frac{partial}{partial x_m} mathcal{E}(x-y)
N(y) dy}), where special integral formulas such as Cauchy formulas
can be used.

	Topics related to this are discussed in Chapter 3 of Part III
of \cite{DS2}, although, as before, further issues are involved there
which are not needed in various settings.  See \cite{BDS} for more
on Clifford analysis, including integral formulas.  Related matters
of curvature are investigated in \cite{H}.

\section{Cauchy integrals and totally real surfaces in ${\bf C}^m$}
\setcounter{equation}{0}

	Let us begin by reviewing some geometrically-oriented linear
algebra, about which Reese Harvey once tutored me.  Fix a positive
integer $m$.  The standard Hermitian inner product on ${\bf C}^m$
is defined by
\begin{equation}
	\langle v, w \rangle = \sum_{j=1}^m v_j \, \overline{w_j},
\end{equation}
where $v$, $w$ are elements of ${\bf C}^m$ and $v_j$, $w_j$ denote their
$j$th components, $1 \le j \le m$.  This expression is complex-linear
in $v$, conjugate-complex-linear in $w$, and satisfies
\begin{equation}
	\langle w, v \rangle = \overline{\langle v, w \rangle}.
\end{equation}
Of course $\langle v, v \rangle$ is the same as $|v|^2$, the square
of the standard Euclidean length of $v$.

	Define $(v, w)$ to be the real part of $\langle v, w \rangle$.
This is a real inner product on ${\bf C}^m$, which is real linear in
both $v$ and $w$, symmetric in $v$ and $w$, and such that $(v, v)$ is
also equal to $|v|^2$.  This is the same as the standard real inner
product on ${\bf C}^m \approx {\bf R}^{2m}$.

	Now define $[v, w]$ to be the imaginary part of $\langle v, w \rangle$.
This is a real linear function in each of $v$ and $w$, and it is
antisymmetric, in the sense that
\begin{equation}
	[w, v] = - [v, w].
\end{equation}
Also, $[v, w]$ is nondegenerate, which means that for each nonzero $v$
in ${\bf C}^m$ there is a $w$ in ${\bf C}^m$ such that $[v, w] \ne 0$.
Indeed, one can take $w = i \, v$.

	Let $L$ be an $m$-dimensional real-linear subspace of ${\bf C}^m$.
We say that $L$ is totally-real if $L$ is transverse to $i \, L$,
where $i \, L = \{ i \, v : v \in L\}$.  Transversality here can be
phrased either in terms of $L \cap i \, L = \{0\}$, or in terms
of $L + i \, L = {\bf C}^m$.

	An extreme version of this occurs when $i \, L$ is the
orthogonal complement of $L$.  Because we are assuming that $L$ has
real dimension $m$, this is the same as saying that elements
of $i \, L$ are orthogonal to elements of $L$.  This is equivalent
to saying that $[v, w] = 0$ for all $v$, $w$ in $L$.  Such a real
$m$-dimensional plane is said to be Lagrangian.

	As a basic example, ${\bf R}^m$ is a Lagrangian subspace
of ${\bf C}^m$.  In fact, the Lagrangian subspaces of ${\bf C}^m$
can be characterized as images of ${\bf R}^m$ under unitary linear
transformations on ${\bf C}^m$.  The images of ${\bf R}^m$ under
special unitary linear transformations, which is to say unitary
transformations with complex determinant equal to $1$, are
called special Lagrangian subspaces of ${\bf C}^m$.

	Now suppose that $M$ is some kind of submanifold or surface in
${\bf C}^m$ with real dimension $m$.  We assume at least that $M$ is a
closed subset of ${\bf C}^m$ which is equipped with a nonnegative
Borel measure $\mu$, in such a way that $M$ is equal to the support of
$\mu$, and the $\mu$-measure of bounded sets are finite.  One might
also ask that $\mu$ behave well in the sense of a doubling condition
on $M$, or even Ahlfors-regularity of dimension $m$.  One may wish to
assume that $M$ is reasonably smooth, and anyway we would ask that $M$
is at least rectifiable, so that $\mu$ can be written as the
restriction of $m$-dimensional Hausdorff measure to $M$ times a
density function, and $M$ has $m$-dimensional approximate tangent
spaces at almost all points.

	Let us focus on the case where $M$ is totally real, so that
its approximate tangent planes are totally real, at least almost
everywhere.  In fact one can consider quantitative versions of this.
Namely, if 
\begin{equation}
	d\nu_m = dz_1 \wedge dz_2 \wedge \cdots \wedge dz_m
\end{equation}
is the standard complex volume form on ${\bf C}^m$, then a linear
subspace $L$ of ${\bf C}^m$ of real dimension $m$ is totally real
if and only if the restriction of $d\nu_m$ to $L$ is nonzero.
In any event, the absolute value of the restriction of $d\nu_m$
to $L$ is equal to a nonnegative real number times the standard
positive element of $m$-dimensional volume on $L$, and positive
lower bounds on that real number correspond to quantitative 
measurements of the extent to which $L$ is totally real.
In the extreme case when $L$ is Lagrangian, this real number
is equal to $1$.  For the surface $M$, one can consider lower
bounds on this real coefficient at each point, or at least 
almost everywhere.

	From now on let us assume that $M$ is oriented, so that the
approximate tangent planes to $M$ are oriented.  This means that
reasonably-nice complex-valued functions on $M$ can be integrated
against the restriction of $d\nu_m$ to $M$.  One can then define
pseudo-accretivity and para-accretivity conditions for the restriction
of $d\nu_m$ to $M$ as in \cite{D-J-S}, which basically mean that
classes of averages of the restriction of $d\nu_m$ to $M$ have nice
lower bounds for their absolute values compared to the corresponding
averages of the absolute value of the restriction to $d\nu_m$ to $m$.
This takes into account the oscillations of the restriction of $d\nu_m$
to $M$.

	Note that if $M$ is a smooth submanifold of ${\bf C}^m$ of
real dimension $m$, then $M$ is said to be Lagrangian if its tangent
spaces are Lagrangian $m$-planes at each point.  This turns out to be
equivalent to saying that $M$ can be represented locally at each point
as the graph of the gradient of a real-valued smooth function on ${\bf
R}^m$ in an appropriate sense, as in \cite{weinstein}.  If the tangent
planes of $M$ are special Lagrangian, then $M$ is said to be a special
Lagrangian submanifold.  See \cite{reese, reese-blaine} in connection
with these.

	It seems to me that there is a fair amount of room here for
various interesting things to come up, basically concerning the
geometry of $M$ and aspects of several complex variables on ${\bf
C}^m$ around $M$.  When $m = 1$, this would include the Cauchy
integral operator applied to functions on a curve and holomorphic
functions on the complement of the curve.  In general this can include
questions about functional calculi, as in \cite{C-M1, C-M2, D-J-S},
and $\overline{\partial}$ problems with data of type $(0,m)$, as well
as relations between the two.

\section{Potentials on various spaces}\footnotetext{A lecture based on this 
section was given at the conference ``Heat kernels and analysis on
manifolds'' at the Institut Henri Poincar\'e, May, 2002.}
\setcounter{equation}{0}

	Let $n$ be a positive integer greater than $1$, and consider
the potential operator $P$ acting on functions on ${\bf R}^n$ defined
by
\begin{equation}
\label{def of P on R^n}
	P(f)(x) = \int_{{\bf R}^n} \frac{1}{|x-z|^{n-1}} \, f(z) \, dz.
\end{equation}
Here $dz$ denotes Lebesgue measure on ${\bf R}^n$.  More precisely,
if $f$ lies in $L^q({\bf R}^n)$, then $P(f)$ is defined almost
everywhere on ${\bf R}^n$ if $1 \le q < n$, it is defined almost
everywhere modulo constants when $q = n$, and it is defined modulo
constants everywhere if $n < q < \infty$.  (If $q = \infty$, then
one can take it to be defined modulo affine functions.)  We shall
review the reasons behind these statements in a moment.

	The case where $n = 1$ is a bit different and special, and we
shall not pay attention to it in these notes for simplicity.
Similarly, we shall normally restrict our attention to functions in
$L^q$ with $1 < q < \infty$.

	A basic fact about this operator on ${\bf R}^n$ is that if $f
\in L^q({\bf R}^n)$, then the first derivatives of $P(f)$, taken in
the sense of distributions, all lie in $L^q({\bf R}^n)$, as long as $1
< q < \infty$.  Indeed, the first derivatives of $P(f)$ are given by
first Riesz transforms of $f$ (modulo normalizing constant factors),
and these are well-known to be bounded on $L^q$ when $1 < q < \infty$.
(In connection with these statements, see \cite{St, SW}.)

	One might rephrase this as saying that $P$ maps $L^q$ into the
Sobolev space of functions on ${\bf R}^n$ whose first derivatives 
lie in $L^q$ when $1 < q < \infty$.  Instead of taking derivatives,
one can look at the oscillations of $P(f)$ more directly, as follows.
Let $r$ be a positive real number, which represents the scale at which
we shall be working.  Consider the expression
\begin{equation}
\label{frac{P(f)(x) - P(f)(y)}{r}}
	\frac{P(f)(x) - P(f)(y)}{r}.
\end{equation}

	To analyze this, let us decompose $P(f)$ into local and
distant parts at the scale of $r$.  Specifically, define operators
$L_r$ and $J_r$ by
\begin{equation}
\label{def of L_r}
	L_r(f)(x) = \int_{\{z \in {\bf R}^n : \, |z-x| < r\}}
			\frac{1}{|x-z|^{n-1}} \, f(z) \, dz
\end{equation}
and
\begin{equation}
\label{def of J_r}
	J_r(f)(x) = \int_{\{z \in {\bf R}^n : \, |z-x| \ge r\}}
			\frac{1}{|x-z|^{n-1}} \, f(z) \, dz.
\end{equation}
Thus $P(f) = L_r(f) + J_r(f)$, at least formally (we shall say more
about this in a moment), so that
\begin{eqnarray}
\label{frac{P(f)(x) - P(f)(y)}{r} = ..., 1}
\lefteqn{\frac{P(f)(x) - P(f)(y)}{r} = }			\\
	&& \frac{L_r(f)(x) - L_r(f)(y)}{r} + \frac{J_r(f)(x) - J_r(f)(y)}{r}.
							\nonumber
\end{eqnarray}

	More precisely, $L_r(f)(x)$ is defined almost everywhere in
$x$ when $f \in L^q({\bf R}^n)$ and $1 \le q \le n$, and it is defined
everywhere when $q > n$.  These are standard results in real analysis
(as in \cite{St}), which can be derived from Fubini's theorem and
H\"older's inequality.  On the other hand, if $1 \le q < n$, then
$J_r(f)(x)$ is defined everywhere on ${\bf R}^n$, because H\"older's
inequality can be used to show that the integral converges.  This does
not work when $q \ge n$, but in this case one can consider the
integral which formally defines the difference $J_r(f)(x) - J_r(y)$.
Namely,
\begin{eqnarray}
\label{J_r(f)(x) - J_r(f)(y), 1}
\lefteqn{\quad J_r(f)(x) - J_r(f)(y) = } 				\\
	& & \int_{{\bf R}^n}
     \biggl(\frac{1}{|x-z|^{n-1}} \, {\bf 1}_{{\bf R}^n \backslash B(x,r)}(z)
  -   \frac{1}{|y-z|^{n-1}} \, {\bf 1}_{{\bf R}^n \backslash B(y,r)}(z) \biggr)
		\, f(z) \, dz.					\nonumber
\end{eqnarray}
Here ${\bf 1}_A(z)$ denotes the characteristic function of a set $A$,
so that it is equal to $1$ when $z \in A$ and to $0$ when $z$ is not
in $A$, and $B(x,r)$ denotes the open ball in ${\bf R}^n$ with center
$x$ and radius $r$.  The integral on the right side of (\ref{J_r(f)(x)
- J_r(f)(y), 1}) does converge when $f \in L^q({\bf R}^n)$ and $q <
\infty$, because the kernel against which $f$ is integrated is bounded
everywhere, and decays at infinity in $z$ like $O(|z|^{-n})$.  This is
easy to check.

	Using this, one gets that $J_r(f)$ is defined ``modulo
constants'' on ${\bf R}^n$ when $f \in L^q({\bf R}^n)$ and $n \le q <
\infty$.  This is also why $P(f)$ can be defined modulo constants on
${\bf R}^n$ in this case (almost everywhere when $q = n$), because of
what we know about $L_r(f)$.  Note that $J_r(f)$ for different values
of $r$ can be related by the obvious formulae, with the differences
given by convergent integrals.  Using this one can see that the
definition of $P(f)$ in terms of $J_r(f)$ and $L_r(f)$ does not depend
on $r$.

	Now let us use (\ref{frac{P(f)(x) - P(f)(y)}{r} = ..., 1})
to estimate $r^{-1} (P(f)(x) - P(f)(y))$.  Specifically, in keeping
with the idea that $P(f)$ should be in the Sobolev space corresponding
to having its first derivatives be in $L^q({\bf R}^n)$ when $f$ is
in $L^q({\bf R}^n)$, $1 < q < \infty$, one would like to see that
\begin{equation}
\label{local average of the difference quotient}
	\frac{1}{|B(x,r)|} \int_{B(x,r)} \frac{|P(f)(x) - P(f)(y)|}{r}
								\, dy
\end{equation}
lies in $L^q({\bf R}^n)$, with the $L^q$ norm bounded uniformly over
$r > 0$.  Here $|A|$ denotes the Lebesgue measure of a set $A$ in
${\bf R}^n$, in this case the ball $B(x,r)$.  In fact, one can even
try to show that the supremum over $r > 0$ of (\ref{local average of
the difference quotient}) lies in $L^q$.  By well-known results, if $q
> 1$, then both conditions follow from the information that the
gradient of $P(f)$ lies in $L^q$ on ${\bf R}^n$, and both conditions
imply that the gradient of $P(f)$ lies in $L^q$.  (Parts of this work
for $q = 1$, and there are related results for the other parts.)  We
would like to look at this more directly, however.

	For the contributions of $L_r(f)$ in (\ref{frac{P(f)(x) -
P(f)(y)}{r} = ..., 1}) to (\ref{local average of the difference
quotient}), one can obtain estimates like the ones just mentioned by
standard means.  For instance, 
\begin{equation}
	\sup_{r > 0} r^{-1} \, L_r(f)(x)
\end{equation}
can be bounded (pointwise) by a constant times the Hardy--Littlewood
maximal function of $f$ (by analyzing it in terms of sums or integrals
of averages of $f$ over balls centered at $x$).  Compare with
\cite{St, SW}.  One also does not need the fact that one has a
difference $L_r(f)(x) - L_r(f)(y)$ in (\ref{frac{P(f)(x) - P(f)(y)}{r}
= ..., 1}), but instead the two terms can be treated independently.
The localization involved is already sufficient to work back to $f$ in
a good way.

	For the $J_r(f)$ terms one should be more careful.  In
particular, it is important that we have a difference $J_r(f)(x) -
J_r(f)(y)$, rather than trying to deal with the two terms separately.
We have seen an aspect of this before, with simply having the
difference be well-defined when $f$ lies in $L^q({\bf R}^n)$ and $n
\le q < \infty$.

	Consider the auxiliary operator $T_r(f)$ defined by
\begin{equation}
\label{def of T_r}
	T_r(f)(x) = \int_{\{z \in {\bf R}^n : \, |z-x| \ge r\}}
			\frac{x-z}{|x-z|^{n+1}} \, f(z) \, dz.
\end{equation}
This is defined everywhere on ${\bf R}^n$ when $f$ lies in $L^q({\bf
R}^n)$ and $1 \le q < \infty$, because of H\"older's inequality.  Note
that $T_r(f)$ takes values in vectors, rather than scalars, because of
the presence of $x-z$ in the numerator in the kernel of the operator.
In fact,
\begin{equation}
    \nabla_x \frac{1}{|x-z|^{n-1}} = - (n-1) \frac{x-z}{|x-z|^{n+1}}.
\end{equation}
Using this and some calculus (along the lines of Taylor's theorem),
one can get that
\begin{eqnarray}
\label{estimate for J_r(f)(x) - J_r(f)(y) - (n-1) (y-x) cdot T_r(f)(x)}
\lefteqn{r^{-1} \, |J_r(f)(x) - J_r(f)(y) - (n-1) (y-x) \cdot
  T_r(f)(x)|}  \\
    & & \qquad\qquad 
 	\le C \int_{{\bf R}^n} \frac{r}{|x-z|^{n+1} + r^{n+1}} 
					\, |f(z)| \, dz	
								\nonumber
\end{eqnarray}
for a suitable constant $C$ and all $x, y \in {\bf R}^n$ with $|x-y|
\le r$.  (In other words, the kernel on the right side of
(\ref{estimate for J_r(f)(x) - J_r(f)(y) - (n-1) (y-x) cdot
T_r(f)(x)}) corresponds to the second derivatives of the kernel of
$J_r$, while $T_r$ reflects the first derivative.)

	The contribution of the right-hand side of (\ref{estimate for
J_r(f)(x) - J_r(f)(y) - (n-1) (y-x) cdot T_r(f)(x)}) to (\ref{local
average of the difference quotient}) satisfies the kind of estimates
that we want, by standard results.  (The right-hand side of
(\ref{estimate for J_r(f)(x) - J_r(f)(y) - (n-1) (y-x) cdot
T_r(f)(x)}) is approximately the same as the Poisson integral of
$|f|$.  Compare with \cite{St, SW} again.)  The remaining piece to
consider is
\begin{equation}
	(n-1) \, r^{-1} \, (y-x) \cdot T_r(f)(x).
\end{equation}
After averaging in $y$ over $B(x,r)$, as in (\ref{local average of the
difference quotient}), we are reduced to looking simply at
$|T_r(f)(x)|$.  Here again the Riesz transforms arise, but in the form
of the truncated singular integral operators, rather than the singular
integral operators themselves (with the limit as $r \to 0$).  By
well-known results, these truncated operators $T_r$ have the property
that they are bounded on $L^q({\bf R}^n)$ when $1 < q < \infty$, with
the operator norm being uniformly bounded in $r$.  Moreover, the
maximal truncated operator
\begin{equation}
	\sup_{r > 0} |T_r(f)(x)|
\end{equation}
is bounded on $L^q({\bf R}^n)$, $1 < q < \infty$.  See \cite{St, SW}.

	These statements are all closely related to the original one
concerning the way that the first derivatives of $P(f)$ are given by
first Riesz transforms of $f$ (up to constant multiples), and lie in
$L^q({\bf R}^n)$ when $f$ does and $1 < q < \infty$.  Instead of
comparing the derivatives of $P(f)$ with Riesz transforms of $f$, we
compare oscillations of $P(f)$ at the scale of $r$ with averages of
$f$ and truncated Riesz transforms of $f$ at the scale of $r$.  We do
this directly, rather than going through derivatives and integrations
of them.

	A nice feature of this discussion is that it lends itself in a
simple manner to more general settings.  In particular, it applies to
situations in which it may not be as convenient to work with
derivatives and integrations of them, while measurements of
oscillations at the scale of $r$ and related estimates still make
sense.

	Instead of ${\bf R}^n$, let us consider a set $E$ in some
${\bf R}^m$.  Let us assume that $E$ is \emph{Ahlfors-regular of
dimension $n$}, by which we mean that $E$ is closed, has at least two
elements (to avoid degeneracies), and that there is a constant $C > 0$
such that
\begin{equation}
\label{Ahlfors-regularity condition}
	C^{-1} \, t^n \le H^n(E \cap \overline{B}(x,t)) \le C \, t^n
\end{equation}
for all $x \in E$ and $t > 0$ with $t \le \diam E$.  Here $H^n$
denotes $n$-dimensional Hausdorff measure (as in \cite{Fe, Ma4}), and
$\overline{B}(x,t)$ denotes the closed ball in the ambient space ${\bf
R}^m$ with center $x$ and radius $t$.

	This condition on $E$ ensures that $E$ behaves
measure-theoretically like ${\bf R}^n$, even if it could be very
different geometrically.  Note that one can have Ahlfors-regular sets
of noninteger dimension, and in fact of any dimension in $(0,m]$ (for
subsets of ${\bf R}^m$).

	Given a function $f$ on $E$, define $P(f)$ on $E$ in the same
manner as before, i.e., by
\begin{equation}
\label{def of P on E}
	P(f)(x) = \int_{E} \frac{1}{|x-z|^{n-1}} \, f(z) \, dz,
\end{equation}
where now $dz$ denotes the restriction of $H^n$-measure to $E$.  Also,
$|x-z|$ uses the ordinary Euclidean distance on ${\bf R}^m$.

	The Ahlfors-regularity of dimension $n$ of $E$ ensures that
$P(f)$ has many of the same basic properties on $E$ as on ${\bf R}^n$.
In particular, if $f$ is in $L^q(E)$, then $P(f)$ is defined almost
everywhere on $E$ (using the measure $H^n$ still) when $1 \le q < n$,
it is defined almost everywhere modulo constants on $E$ when $q = n$,
and it is defined everywhere on $E$ modulo constants when $n < q <
\infty$.  One can show these statements in essentially the same manner
as on ${\bf R}^n$, and related results about integrability, bounded
mean oscillation, and H\"older continuity can also be proven in
essentially the same manner as on ${\bf R}^n$.  

	What about the kind of properties discussed before, connected
to Sobolev spaces?  For this again one encounters operators on
functions on $E$ with kernels of the form
\begin{equation}
\label{kernels for operators}
	\frac{x-z}{|x-z|^{n+1}}.
\end{equation}
It is not true that operators like these have the same kind of
$L^q$-boundedness properties as the Riesz transforms do for arbitrary
Ahlfors-regular sets in ${\bf R}^m$, but this is true for integer
dimensions $n$ and ``uniformly rectifiable'' sets $E$.  In this
connection, see \cite{Ca, CDM, CMM, Da1, Da2, Da3, Da4, DS1, DS2, Ma4,
MMV}, for instance (and further references therein).

	When $E$ is not a plane, the operators related to the kernels
(\ref{kernels for operators}) are no longer convolution operators, and
one loses some of the special structure connected to that.  However,
many real-variable methods still apply, or can be made to work.  See
\cite{CW1, CW2, CM*, jean-lin}.  For example, the Hardy--Littlewood maximal
operator still behaves in essentially the same manner as on Euclidean
spaces, as do various averaging operators (as were used in the earlier
discussion).  Although one does not know that singular integral
operators with kernels as in (\ref{kernels for operators}) are bounded
on $L^q$ spaces for arbitrary Ahlfors-regular sets $E$, there are
results which say that boundedness on one $L^q$ space implies
boundedness on all others, $1 < q < \infty$.  Boundedness of singular
integral operators (of the general Calder\'on--Zygmund type) implies
uniform boundedness of the corresponding truncated integral operators,
and also boundedness of the maximal truncated integral operators.

	At any rate, a basic statement now is the following.  Let $n$
be a positive integer, and suppose that $E$ is an Ahlfors-regular set
in some ${\bf R}^m$ which is ``uniformly rectifiable''.  Define the
potential operator $P$ on functions on $E$ as in (\ref{def of P on
E}).  Then $P$ takes functions in $L^q(E)$, $1 < q < \infty$, to
functions on $E$ (perhaps modulo constants) which satisfy ``Sobolev
space'' conditions like the ones on ${\bf R}^n$ for functions with
gradient in $L^q$.  In particular, one can look at this in terms of
$L^q$ estimates for the analogue of (\ref{local average of the
difference quotient}) on $E$, just as before.  These estimates can be
derived from the same kinds of computations as before, with averaging
operators and operators like $T_r$ in (\ref{def of T_r}), but now on
$E$.  The estimates for $T_r$ use the assumption of uniform
rectifiability of $E$ (boundedness of singular integral operators).
The various other integral operators, with the absolute values inside
the integral sign, are handled using only the Ahlfors-regularity of
$E$.

	Note that for sets $E$ of this type, one does not necessarily
have the same kind of properties concerning integrating derivatives as
on ${\bf R}^n$.  In other words, one does not automatically get as
much from looking at infinitesimal oscillations, along the lines of
derivatives, as one would on ${\bf R}^n$.  The set $E$ could be quite
disconnected, for instance.  However, one gets the same kind of
estimates at larger scales for the potentials that one would normally
have on ${\bf R}^n$ for a function with its first derivatives in
$L^q$, by looking at a given scale $r$ directly (rather than trying to
integrate bounds for infinitesimal oscillations), as above.

	For some topics related to Sobolev-type classes on general
spaces, see \cite{FHK, Ha, HaK1, HaK2, HeK1, HeK2} (and references
therein).

	Although the potential operator in (\ref{def of P on E}) has a
nice form, it is also more complicated than necessary.  Suppose that
$E$ is an $n$-dimensional Lipschitz graph, or that $E$ is simply
bilipschitz--equivalent to ${\bf R}^n$, or to a subset of ${\bf R}^n$.
In these cases the basic subtleties for singular integral operator
with kernel as in (\ref{kernels for operators}) already occur.
However, one can obtain potential operators with the same kind of nice
properties by making a bilipschitz change of variables into ${\bf
R}^n$, and using the classical potential operator there.  This leads
back to the classical first Riesz transforms on ${\bf R}^n$, as in
\cite{St, SW}.

	Now let us consider a rather different kind of situation.
Suppose that $E$ is an Ahlfors-regular subset of dimension $n$ of some
${\bf R}^m$ again.  For this there will be no need to have particular
attention to integer values of $n$.  Let us say that $E$ is a
\emph{snowflake} of order $\alpha$, $0 < \alpha < 1$, if there is a
constant $C_1$ and a metric $\rho(x,y)$ on $E$ such that
\begin{equation}
\label{snowflake condition}
	C_1^{-1} \, |x-y| \le \rho(x,y)^\alpha \le C_1 \, |x-y|
\end{equation}
for all $x, y \in E$.

	In this case, let us define a potential operator
$\widetilde{P}$ on functions on $E$ by
\begin{equation}
\label{def of widetilde{P}}
	\widetilde{P}(f)(x) 
	  = \int_E \frac{1}{\rho(x,z)^{\alpha (n-1)}} \, f(z) \, dz.
\end{equation}
Here $dz$ denotes the restriction of $n$-dimensional Hausdorff measure
to $E$ again.  This operator is very similar to the one before, since
$\rho(x,z)^{\alpha (n-1)}$ is bounded from above and below by constant
multiples of $|x-z|^{n-1}$, so that the kernel of $\widetilde{P}$ is
bounded from above and below by constant multiples of the kernel of
the operator $P$ in (\ref{def of P on E}).

	This operator enjoys the same basic properties as before, with
$\widetilde{P}(f)$ being defined almost everywhere when $f$ lies in
$L^q(E)$ and $1 \le q < n$, defined modulo constants almost everywhere
when $q = n$, and defined modulo constants everywhere when $n < q <
\infty$, for essentially the same reasons as in the previous
circumstances.  However, there is a significant difference with this
operator, which one can see as follows.  Let $x$, $y$, $z$ be three
points in $E$, with $x \ne z$ and $y \ne z$.  Then
\begin{equation}
\label{inequality for differences of kernels}
 	\biggl| \frac{1}{\rho(x,z)^{\alpha (n-1)}}  
                  - \frac{1}{\rho(y,z)^{\alpha (n-1)}} \biggr|
 \le C \, \frac{\rho(x,y)}{\min(\rho(x,z),\rho(y,z))^{\alpha (n-1) + 1}}
\end{equation}
for some constant $C$ which does not depend on $x$, $y$, or $z$,
but only on $\alpha (n-1)$.  Indeed, one can choose $C$ so that
\begin{equation}
\label{inequality simply for positive numbers}
	\bigl| a^{\alpha (n-1)} - b^{\alpha (n-1)} \bigr|
		\le C \, \frac{|a-b|}{\min(a,b)^{\alpha (n-1) + 1}}
\end{equation}
whenever $a$ and $b$ are positive real numbers.  This is an elementary
observation, and in fact one can take $C = \alpha (n-1)$.  One can get
(\ref{inequality for differences of kernels}) from (\ref{inequality
simply for positive numbers}) by taking $a = \rho(x,z)$ and $b =
\rho(y,z)$, and using the fact that
\begin{equation}
	| \rho(x,z) - \rho(y,z)| \le \rho(x,y).
\end{equation}
This last comes from the triangle inequality for $\rho(\cdot, \cdot)$,
which we assumed to be a metric.

	Using the snowflake condition (\ref{snowflake condition}), we
can obtain from (\ref{inequality for differences of kernels}) that
\begin{equation}
\label{inequality for differences of kernels, 2}
 	\biggl| \frac{1}{\rho(x,z)^{\alpha (n-1)}}  
                  - \frac{1}{\rho(y,z)^{\alpha (n-1)}} \biggr|
  \le C' \, \frac{|x-y|^{1/\alpha}}{\min(|x-z|, |y-z|)^{(n-1) + 1/\alpha}}
\end{equation}
for all $x, y, z \in {\bf R}^n$ with $x \ne z$, $y \ne z$, and with a
modestly different constant $C'$.  The main point here is that the
exponent in the denominator on the right side of the inequality is
strictly larger than $n$, because $\alpha$ is required to lie in
$(0,1)$.  In the previous contexts, using the kernel $1/ |x-z|^{n-1}$
for the potential operator, there was an analogous inequality with
$\alpha = 1$, so that the exponent in the denominator was equal to
$n$.

	With an exponent larger than $n$, there is no need for
anything like singular integral operators here.  More precisely, there
is no need for the operators $T_r$ in (\ref{def of T_r}) here; one can
simply drop them, and estimate the analogue of $|J_r(f)(x) -
J_r(f)(y)|$ when $|x - y| \le r$ directly, using (\ref{inequality for
differences of kernels, 2}).  In other words, one automatically gets
an estimate like (\ref{estimate for J_r(f)(x) - J_r(f)(y) - (n-1)
(y-x) cdot T_r(f)(x)}) in this setting, without the $T_r$ term, and
with some minor adjustments to the right-hand side.  Specifically, the
$r$ in the numerator on the right side of (\ref{estimate for J_r(f)(x)
- J_r(f)(y) - (n-1) (y-x) cdot T_r(f)(x)}) would become an
$r^{1/\alpha -1}$ in the present situation, and the exponent $n+1$ in
the denominator would be replaced with $n - 1 + 1/\alpha$.  This leads
to the same kinds of results in terms of $L^q$ norms and the like as
before, because the rate of decay is enough so that the quantities in
question still look like suitable averaging operators in $f$.  (That
is, they are like Poisson integrals, but with somewhat less decay.
The decay is better than $1/|x-z|^n$, which is the key.  As usual, see
\cite{St, SW} for similar matters.)

	The bottom line is that if we use the potential operator
$\widetilde{P}$ from (\ref{def of widetilde{P}}) instead of the
operator $P$ from (\ref{def of P on E}), then the two operators are
approximately the same in some respects, with the kernels being of
comparable size in particular, but in this situation the operator
$\widetilde{P}$ has the nice feature that it automatically enjoys the
same kind of properties as in the ${\bf R}^n$ case, in terms of
estimates for expressions like (\ref{local average of the difference
quotient}) (under the snowflake assumption for $E$).  That is, one
automatically has that $\widetilde{P}(f)$ behaves like a function in a
Sobolev class corresponding to first derivatives being in $L^q$ when
$f$ lies in $L^q$.  One does not need $L^q$ estimates for singular
integral operators for this, as would arise if we did try to use the
operator $P(f)$ from (\ref{def of P on E}).

	These remarks suggest numerous questions...

	Of course, some other basic examples involve nilpotent Lie
groups, like the Heisenberg group, and their invariant geometries.

	As a last comment, note that for the case of snowflakes we
never really needed to assume that $E$ was a subset of some ${\bf
R}^m$.  One could have worked just as well with abstract metric spaces
(still with the snowflake condition).  However, Assouad's embedding
theorem \cite{A1, A2, A3} provides a way to go back into some ${\bf
R}^m$ anyway.  The notion of uniform rectifiability makes sense for
abstract metric spaces, and not just subsets of ${\bf R}^m$, and an
embedding into some ${\bf R}^m$ is sometimes convenient.  In this
regard, see \cite{S5}.

\end{document}